\documentclass{amsart}
\usepackage[latin1]{inputenc}
\usepackage{amssymb}
\usepackage{latexsym}
\usepackage{longtable}

\newtheorem{teor}{Theorem}[section]

\newtheorem{lema}[teor]{Lemma}
\newtheorem{prop}[teor]{Proposition}

\newtheorem{remark}{Remark}[section]

\newcommand{\Q}{\mathbb{Q}}
\newcommand{\Z}{\mathbb{Z}}
\newcommand{\F}{\mathbb{F}}

\newcommand{\C}{\mathbb{C}}

\newcommand{\GL}{\operatorname{GL}}

\newcommand{\End}{\operatorname{End}}

\newcommand{\Jac}{\operatorname{Jac}}

\newcommand{\id}{\operatorname{id}}
\newcommand{\Aut}{\operatorname{Aut}}

\newfont{\gotip}{eufb10 at 12pt}

\newcommand{\om}{{\omega }}
\newcommand{\ra}{{\rightarrow }}

\newcommand{\R}{{\mathbb R}}

\include{thebibliography}

\begin{document}

\title[Equations of Shimura curves of genus $2$]
{Equations of Shimura curves of genus $2$}

\author{Josep Gonz\'{a}lez, Victor Rotger}
\footnote{The first author is supported in part by   DGICYT Grant
BFM2000-0794-C02-02 and the second author is partially supported by Ministerio
de Ciencia y Tecnolog\'{\i}a BFM2000-0627}

\address{Universitat Polit\`{e}cnica de Catalunya,
Departament de Matem\`{a}tica Aplicada IV (EUPVG), Av.\ Victor
Balaguer s/n, 08800 Vilanova i la Geltr\'{u}, Spain.}

\email{josepg@mat.upc.es,vrotger@mat.upc.es}



\begin{abstract}
We present explicit models for Shimura curves $X_D$ and
Atkin-Lehner quotients $X_D/\langle \om _m \rangle $ of them of
genus $2$. We show that several equations conjectured by Kurihara
are correct and compute for them the kernel of Ribet's isogeny
$J_0(D)^{new}\rightarrow J_D$ between the new part of the Jacobian
of the modular curve $X_0(D)$ and the Jacobian of $X_D$.
\end{abstract}

\maketitle

\section{Introduction}

Let $B_D$ be the indefinite quaternion algebra over $\Q $ of
reduced discriminant $D=p_1\cdot ... \cdot p_{2r}$ for pairwise
different prime numbers $p_i$ and let $X_D/\Q $ be the Shimura
curve attached to $B_D$. As it was shown by Shimura \cite{Sh67},
$X_D$ is the coarse moduli space of abelian surfaces with
quaternionic multiplication by $B_D$.

Let $W=\{ \om _m: m\mid D\}\subseteq \Aut _{\Q }(X_D)$ be the
group of Atkin-Lehner involutions. For any $m\mid D$, we shall
denote $X_D^{(m)}=X_D/\langle \om _m\rangle $ the quotient of the
Shimura curve $X_D$ by $\om _m$. The importance of the curves
$X_D^{(m)}$ is enhanced by their moduli interpretation as curves
embedded in Hilbert-Blumenthal surfaces and Igusa's threefold
$\mathcal A _2$ (cf.\,\cite{Ro2},\,\cite{Ro3}).

The classical modular case arises when $D=1$. In this case,
automorphic forms of these curves admit Fourier expansions around
the cusp of infinity and we know explicit generators of the field
of functions of such curves. Also, explicit methods are known to
determine bases of the space of regular differentials of them,
which are used to compute equations for quotients of modular
curves.

When $D\not =1$, the absence of cusps has been an obstacle for
explicit approaches to Shimura curves. Explicit methods to handle
with functions and regular differential forms on these curves are
less accessible and we refer the reader to \cite{Bay} for progress
in this regard. For this reason, at present few equations of
Shimura curves are known, all of them of genus $0$ or $1$ (cf.
 \cite{Kur}, \cite{JoLi}, \cite{Elk}). In addition, in a later work, Kurihara
conjectured equations for all Shimura curves of genus $2$ and for
several curves of genus $3$ and $5$, though he was not able to
give a proof for his guesses (cf. \cite{Kur2}).

In this paper, we present equations for thirteen genus two
bielliptic Shimura curves and Atkin-Lehner quotients of them.

In particular, we prove that the equations suggested in
\cite{Kur2} for $X_{26}$, $X_{38}$ and $X_{58}$ are
unconditionally correct.
In turn, this has allowed us to explicitly determine the kernel of
Ribet's isogeny $J_0(D)^{new}\ra \Jac (X_D)$ and to prove that
Ogg's prediction in \cite{Ogg2} is also correct for these cases.

The remaining $10$ curves presented here are the only bielliptic
curves $X_D^{(m)}$, $m\not =1$, of genus $2$. Phrased in other
terms, this is the complete list of all genus two curves
$X_D^{(m)}$ whose hyperelliptic involution is not of Atkin-Lehner
type. Note that a phenomenon of this kind was already encountered
in the modular setting by the curve $X_0(37)$.
Our method can also be used to determine equations for genus two
bielliptic Shimura curves with nontrivial $\Gamma _0$-level
structure, of which there exist $89$. For the sake of
brevity, these will not be considered in this work.

\section{Explicit models of bielliptic curves of genus $2$}

\begin{prop}\label{explicit}
Let $C$ be a genus two curve defined over a field $k$ of
characteristic not $2$ or $3$ and $w$ its hyperelliptic
involution. Assume $\Aut_k (C)$ contains a subgroup $\langle
u_1,u_2=u_1. w\rangle $ isomorphic to $(\Z/2\Z)^2$ and let us
denote by $C_{u_i}$ the elliptic quotient $C/\langle u_i\rangle$.
If the two elliptic curves
$$
E_1: Y^2=X^3+ A_1 X+B_1\,,\quad E_2: Y^2=X^3+ A_2 X+B_2\,,
$$
are isomorphic to $C_{u_1}$ and $C_{u_2}$ resp. over $k$, then $C$
admits a hyperelliptic equation  of the form  $y^2=a x^6+ b x^4+ c
x^2 +d$, where $a\in k^*$, $b\in k$ are solutions of the following
system:
\begin{equation}\label{eq1}
\left \{
\begin{array}{ccc}
27 a^3 B_2&= &2A_1^3+ 27 B_1^2 +  9A_1 B_1 b+ 2A_1^2b^2 -B_1 b^3
\\  9 a^2 A_2&=&-3 A_1^2 + 9 B_1 b+ A_1 b^2
\end{array} \right. \,,
\end{equation}
$c=(3A_1 + b^2)/(3a)$, $ d=(27 B_1+9 A_1 b + b^3)/(27 a^2)$ and
the involution $u_1$ on $C$ is given by $(x,y)\mapsto (-x,y)$.

\end{prop}

\noindent

{\bf Proof.} If  $C$ has a nonhyperelliptic
involution $u_1$ defined over $k$, then $C/k$ admits  an equation
of the form
\begin{equation}\label{eq2}
y^2 = a x^6+ b x^4+ c x^2 + d\,, \quad a,b,c,d\in k\,,
\end{equation}
and the involution $u_1$ acts sending  $(x,y)\mapsto (-x,y)$.
Indeed, due to the fact that the morphisms $C\rightarrow C/\langle u_i\rangle$ are defined over $k$ there are $\omega_1, \omega_2\in
H^0(C,\Omega_{C/k}^1)$ such that  $u_1^{*}\omega_1=\omega_1$,
$u_1^{*}\omega_2=-\omega_2$.
The functions
$x=\omega_1/\omega_2$, $y= dx/\omega_2$ must satisfy a relation of
the form  $y^2=f(x)$, where $f\in k[x]$ has degree $5$ or $6$ and
does not have double roots (see Proposition 2.1 in \cite{gogo}). Then, $u_1^{*}(x)=-x$ and
$u_1^{*}(y)=y$. It follows $f(-x)=x$ and, in particular, $\deg
f=6$.

Given an equation for $C$ as (\ref{eq2}), the elliptic curves
$$
C_1: Y^2= a X^3+ b X^2+ c X+d\,,\quad C_2: Y^2= d X^3+ c X^2+ b X+
a
$$
are $k$-isomorphic to $C_{u_1}$ i $C_{u_2}$ respectively, due to
the fact that the nonconstant morphisms $\pi_1: C\rightarrow C_1$,
$(x,y)\mapsto (x^2\,, \,y)$ i $\pi_2: C\rightarrow C_2$,
$(x,y)\mapsto (1/x^2\,, \,\, y/x^3)$ are defined over $k$ and
satisfy  $\pi_i\circ u_i=\pi_i$, $i\leq 2$. Therefore, every curve
$C_{u_i}$ is isomorphic over $k$ to the curve $ Y^2= X^3+ A_{u_i}
X+ B_{u_i}$, where
\begin{equation}\label{sist}\left \{
\begin{array}{cc}
 A_{u_1}=-b^2/3 + a.c\,,  & B_{u_1}=2\,b^3/27 - a.b.c/3 + a^2.d\,,\\
  A_{u_2}=-c^2/3 + b.d \,, & B_{u_2}=2\,c^3/27 - b.c.d/3 + a.d^2\,,
  \end{array}
  \right.
\end{equation}
and, thus, there exist $\mu_i \in k^{*}$, $i\leq 2$, such that:
$$
A_{u_i} \mu_i^4=A_i\,,\quad B_{u_i} \mu_i^6=B_i \,.
$$
It can be easily checked that the curve
$$
y^2= a \frac{\mu_1^4}{\mu_2^2} x^6+ b\mu_1^2 x^4+ c \mu_2 ^2 x^2
+d \frac{\mu_2 ^4}{\mu_1^2}
$$
is $k$-isomorphic to $C$. The statement is an immediate
consequence of rewriting the system (\ref{sist}) for this
equation, since $a\neq 0$ and now
   $  A_{u_i} =A_i$, $B_{u_i} =B_i$.  $\Box$

\begin{remark}
Given two elliptic curves $E_1$, $E_2$ over $k$ and a group
isomorphism $\psi:E_1[2](\bar k)\rightarrow E_2[2](\bar k)$ which
is not the restriction of an isomorphism between $E_1$ and $E_2$
over $\bar k$, Proposition 4 in \cite{HoLePo} yields a genus two
curve $C/\bar k$ such that $\Jac C\simeq (E_1\times E_2
)/\{(P,\psi (P)):P\in E_1[2]\}$.

In our case, when we consider the elliptic curves defined over $k$
$$
C_1: Y^2= a X^3+ b X^2+ c X+d\,,\quad C_2: Y^2= d X^3+ c X^2+ b X+
a
$$
and the isomorphism of $G_k$-modules $\psi: C_1[2](\bar
k)\rightarrow C_2[2](\bar k)$, $(x,0)\mapsto (1/x,0)$, the formula
of the quoted proposition yields a curve which is shown to be
isomorphic to $C:y^2=a x^6+ b x^4+c x^2+d$ over $k$.

Hence, system (\ref{eq1}) can be viewed as a different way to
collect all curves $C$ obtained from all $\psi$ as above. By
Proposition 3 of \cite{HoLePo}, if $E_1\not \simeq E_2$ over $\bar
k$, the system has six different solutions $(a,b,c,d)\in\bar k^4$
and there is a unique solution defined over $k$ if and only
$(E_1\times E_2)(\bar k) $ has a unique nontrivial $G_k$-stable
subgroup $G$ isomorphic to $(\Z/2\Z )^2$, and in this case $\Jac
C=(E_1\times E_2)/G$.  Here, by the trivial $G_k$-stable
subgroups, we mean $E _1[2](\bar k)\times \{0\}$ and $\{0\}\times
E _2[2](\bar k)$.
\end{remark}

\section{Shimura curves of genus two}

Let $B$ be an indefinite quaternion algebra over $\Q $ of
discriminant $D$ and let $X_D$ denote the Shimura curve attached
to it. For any integer $N\geq 1$, let $X_0(N)$ be the modular
curve of level $N$ and $J_0(N) = \Jac (X_0(N))$. By
$J_0(N)^{new}$, we shall denote the new part of $J_0(N)$ viewed as
an optimal quotient of it. Ribet's isogeny theorem establishes the
existence of a Hecke invariant isogeny
$$J_0(D)^{new}\ra \Jac (X_D)$$ over $\Q $, though its proof relies on
the fact that both abelian varieties have the same $L$-series and
therefore is not explicit (cf.\,\cite{Ri}, see also \cite{Ar}).

The problem of determining the possible kernels of the isogeny has
been studied by Ogg in \cite{Ogg2}, the underlying idea being that
the knowledge of the group of connected components of the N\'{e}ron
models of $J_0(D)^{new}$ and $\Jac (X_D)$ at a prime $p\mid D$
yields necessary conditions to be satisfied by the isogenies
between them. As in \cite{Ogg2}, the component groups of $\Jac
(X_D)$ can be handled by Raynaud's method and the theory of
\v{C}erednik-Drinfeld. However, the component groups of the optimal
quotients $J_0(D)^{new}$ were only recently determined by Conrad
and Stein in \cite{CoSt}.

The aim of this section is to provide equations for the three
Shimura genus two curves and to make Ribet's isogeny explicit for
these examples.

\begin{teor}\label{26}

The curves $X_D$ with $D=26,38,58$ are the unique Shimura curves
of genus two. Moreover,

\begin{enumerate}

\item[(i)]
Equations for the curves $X_D$ are given in the following table:
$$
\begin{array}{|c|c|}
\hline
D & X_D \\
\hline
26 & y^2 =-2 x^6+19 x^4- 24 x^2-169  \\
\hline
38 & y^2=-19 x^6 -82 x^4 -59 x^2 -16  \\
\hline
58 & 2 y^2=  - x^6-39 x^4- 431x^2- 841  \\
\hline

\end{array}
$$

\item[(ii)]
In all of three cases $J_0(D)^{new}$ is the Jacobian of a genus
two curve $C_D$ defined over $\Q$ and there is a cuspidal divisor
$c(D)$ in $J_0(D)^{new}$ and an exact sequence
$$0\ra \langle c(D)\rangle \ra J_0(D)^{new}\ra \Jac (X_D)\ra 0.$$

Equations for the curves $C_D$, the cuspidal divisors $c(D)$ and
their orders are given in the following table:

$$
\begin{array}{|c|c|c|c|}
\hline
D & C_D & c(D) & |\langle c(D)\rangle |\\
\hline
26 &  y^2=   13x^6  + 10x^4 - 3x^2-4& (1/13)-(\infty) & 7\\
\hline
38 &  y^2 =   x^6  + 2x^4  + x^2+ 76 & (1/19)-(\infty) & 5 \\
\hline
58 &   y^2=     x^6 + 6x^4- 7x^2+16& (1/29)-(\infty) & 5\\
\hline
\end{array}
$$

\end{enumerate}

\end{teor}

\noindent {\bf Proof.} It follows from Ogg's list of hyperelliptic
Shimura curves (cf. \cite{Ogg1}) that $D=26, 38$ and $58$ are the
only values of $D$ for which $g(X_D)=2$. These curves are
bielliptic; more precisely, in Cremona's notation, by \cite{Rob}, it follows that for
these values of $D$, $X_D/\langle w_2\rangle$ is the elliptic
curve $B2$ of conductor $D$ while $X_{26}/\langle w_{13}\rangle$,
$X_{38}/\langle w_{19}\rangle$, $X_{58}/\langle w_{58}\rangle$ are
the elliptic curves $26A_1$, $38A_1$ and $58A_1$, respectively. It
can be checked that for these values of $D$, the classes of
isomorphism over $\bar \Q$ of both curves are different. Applying
Proposition \ref{explicit}, we obtain that in all these cases the
system (\ref{eq1}) gives a unique genus two curve defined over
$\Q$, which is given in the first table of the statement.

Let $f_1$ and $f_2$ be the two normalized newforms of
$S_2(\Gamma_0(D))$ and let  $E_A$ and $E_B$ be the elliptic curves
over $\Q$ which are the strong Weil curve in the class of isogeny
$A$ and $B$ respectively. We know that $J_0(D)^{new}$ and
$E_A\times E_B$ are isogenous over $\Q$. We compute the following
lattices of $\C$:
$$
\Lambda_i=\{ \int_{\gamma}f_i(q)dq/q:\gamma \in H_1(X_0(D),\Z)\}\,,
\quad 1\leq i\leq 2\,,
$$
and the lattice of $\C ^2$:
$$
\Lambda=\{
(\int_{\gamma}f_1(q)dq/q,\int_{\gamma}f_2(q)dq/q):\gamma \in
H_1(X_0(D),\Z)\}.
$$
We obtain $(\Lambda_1\times \Lambda_2)/\Lambda\simeq (\Z/2\Z)^2$
with $\Lambda$ being  different from the lattices $ (1/2\cdot
\Lambda_1)\times \Lambda_2$ and $\Lambda_1\times(1/2\cdot
\Lambda_2)$. This result implies that there exists a nontrivial
$G_{\Q}$-stable subgroup $G$ of $E_A\times E_B$ isomorphic to
$(\Z/2\Z)^2$   such that $J_0(D)^{new}=(E_A\times E_B)/G$. In
consequence, by \cite{HoLePo}, we know that $J_0(D)^{new}$ is the
Jacobian of a genus two curve $C_D/\Q$ (of course, for $D=26$, it
was already known). Again, applying Proposition \ref{explicit}, we
obtain  a unique genus two curve defined over $\Q$, which is given
in the table of the statement.  Note that equations for $X_0(26)$
were already known (cf.\ \cite{Go}).

Now, we consider the morphism $\phi$ obtained as the composition
of the following morphisms defined over $\Q$
$$
J_0(D)^{new}\stackrel{\mu}\longrightarrow E_A\times
E_B\stackrel{\id_A\times \phi_B}\longrightarrow E_{A}\times
E_{B_2}\stackrel{\nu}\longrightarrow \Jac (X_D)\,,
$$
where $\ker\mu,\ker \nu\simeq (\Z/2\Z)^2$, $\id_A$ is the identity on $E_A$ and $\phi_B$ is
the cyclic isogeny from $E_B$ to $E_{B_2}$. One  can check (see
\cite{Cre}) that in all these cases the group $\ker
(\id_A\times\phi_B )$ is a subgroup of $E_A(\Q)\times E_B(\Q)$ of
cardinality $ 7$, $ 5$, $ 5$ depending on
whether $D$ is $26$, $38$ or $58$  and, moreover, this group is
the unique subgroup of rational points of $E_A(\Q)\times E_B(\Q)$
with such a cardinality. Since $\id_A\times \phi_B $ has odd
degree, this morphism  maps the kernel of $\hat \mu$ to the kernel
of $\nu$ because both kernels are the unique nontrivial
$G_{\Q}$-stable subgroups isomorphic to $(\Z/2\Z)^2$ in their
abelian varieties. Then, there is a morphism
$\phi':J_0(D)^{new}\rightarrow \Jac (X_D)$ such that $
\phi=[2]\phi'$ and, thus, $|\ker \phi'|=|\ker (\id_A\times
\phi_B)|=|\ker (\phi_B)|$.

Since $D$ is square-free, we recall that cuspidal divisors are
rational points in $J_0(D)$ and in particular in $J_0(D)^{new}$.
The cuspidal divisors $c(D)$ given in the table have order $7$,
$15$ and $35$ in $J_0(D)$ for $D=26,38$ and $58$ respectively. For
the cases $D=38$ and $58$, we compute $(\int_{c(D)}f_1(q)
dq/q,\int_{c(D)}f_2(q) dq/q)\in \Q\otimes \Lambda$ and  check that
its order in $(\Q\otimes \Lambda)/\Lambda$  is $5$. This concludes
the proof. $\Box$

\begin{remark}
The three equations agree with the equations suggested in
\cite{Kur2}. Moreover, Ogg suggested in \cite[p.\ 213]{Ogg2}, that
the minimal degree of Ribet's isogeny should be the numerator of
$\frac{p+1}{12}$ whenever $D=2 p$. This agrees with the table
above.
More precisely, when $D=26$,
$(0)-(1/2)+(1/13)-(\infty)=2(1/13)-2(\infty)$ in $J_0(26)$, which
proves that the prediction done by Ogg in \cite{Ogg2} about the kernel of this isogeny is again
right.

It can also be checked that for $D=26,38$ the group
$\langle c(D)\rangle $ is generated by $3(0)-3(\infty)$ while for
$D=58$ it is generated by $(0)-(\infty)$, and in all three cases
the kernel of the isogeny is a subgroup of $\langle
(0)-(\infty)\rangle$. It would be interesting to know whether the
pattern suggested by the examples holds in greater generality.
\end{remark}

\begin{remark}
Theorem \ref{26} provides an explicit model $y^2=ax^6+bx^4+cx^2+d$
for $X_{26}$, $X_{38}$ and $X_{58}$ which is known to have a cusp
singuarity at the only point $P_{\infty }=[0:1:0]$ of infinity. A
smooth model of the curve is obtained by blowing up the point; the
preimage of $P_{\infty }$ by the normalizing map are two points and the coordinates of everyone of them
generates $\Q (\sqrt {a})$. In the three cases
above $\Q (\sqrt {a})$ is quadratic imaginary, as it was expected
since Shimura curves fail to have real points \cite{Sh2}.
\end{remark}

\section{Explicit models of Atkin-Lehner quotients of
Shimura curves}

Let $D=p_1\cdot ...\cdot p_{2 r}$, $r\geq 1$, and $m\mid D$. Let
$X_D^{(m)} = X_D/\langle \om _m\rangle $ be the quotient of the
Shimura curve $X_D$ by the Atkin-Lehner involution $\om _m$.

Let $\mathbb T = \langle T_{\ell }, \om _p: \ell \nmid D, p\mid D
\rangle _{\Q }$ and $\underline{\mathbb T} = \langle
\underline{T}_{\ell }, \underline{\om }_p: \ell \nmid D, p\mid D
\rangle _{\Q }$ denote the Hecke algebra regarded as $\mathbb T =
\End _{\Q }(\Jac (X_D))\otimes \Q $ and $\underline{\mathbb T}
=\End _{\Q }(J_0(D)^{new})\otimes \Q $, respectively. Ribet's
isogeny $\Jac (X_D)\ra J_0(D)^{new}$ provides an isomorphism between the vector spaces of regular differentials and identifies
$T_{\ell }$ with $\underline{T}_{\ell }$ and $\om _m$ with $\mu
(m) \underline{\om }_m$ for any $m\mid D$, where $\mu
(m)=(-1)^{\sharp \{ \mathrm{primes }\ p\mid m \}}$.

\begin{lema}\label{g2}
The genus of $X_D^{(m)}$ is $2$ if and only if $(D,m)\in \{
(35,5),(39,3),(51,17),\\(55,11),(57,3),
(62,31),(65,5),\,(65,13),(69, 23),(74,2),(74,37),(82,2),\,(85, 5),
(85,\\ 85),\, (86, 2),\,(86,43),\,(87,3),\,(91,91),\,
(93,93),\,(94,47),\,(106,2),\,(115,115),\,(118,2),\\(122,61),\,(123,123),
(129,43),\,
(141,141),\,(142,2),\,(142,142),\,(155,155),\,(158,158),\\(161,
161),(166,83),(178,178),(183,183),(237,79),(254,254),(326,326),(446,446)\}.$
\end{lema}
\noindent {\bf Proof.} Assume that the pair $(D, m)$ is such that
$g(X_D^{(m)})=2$. Since $\Aut (X_D^{m})\supseteq W/\langle \om _m
\rangle \simeq (\Z /2\Z )^{2 r -1}$ and curves of genus two
contain at most two copies of the cyclic group of order $2$, it
follows that necessarily $r=1$ and hence $D=p\cdot q$.

Let $\ell \nmid D$ be a prime of good reduction of the  curve.
Following \cite[\S5]{Ogg1}, we obtain $\varphi (D)(\ell
-1)/12\leq |\widetilde{X}_D (\F_{\ell ^2})|\leq 2
|\widetilde{X}_D^{(m)} (\F_{\ell ^2})|$, where $\widetilde{X}_D$
denotes the special fiber of Morita's integral model of $X_D$ over
$\Z _{\ell }$. Since $\widetilde{X}_D^{(m)}$ is hyperelliptic, it
admits a map of degree $2$ onto the projective line and hence
$|\widetilde{X}_D^{(m)}(\F _{\ell ^2})|\leq 2 (\ell ^2 + 1)$.
Thus, $\varphi (D)\leq 48 (\ell ^2 + 1)/(\ell -1)$. Since
$g(X_6)=0$, we may choose either $\ell =2$ or $\ell = 3$ and hence
$\varphi (D)\leq 240$. A computation of genera now yields
the lemma. $\Box $

\begin{prop}\label{AL}

A Shimura curve $X_D^{(m)}$ of genus two admits a bielliptic
involution if and only if $(D, m) \in \{(91, 91), (123, 123),\,
(141, 141),\, (142, 2),\, (142, 142),\,(155, \\ 155),(158, 158),
(254, 254), (326, 326), (446,446)\} $.

In all these $10$ cases, the hyperelliptic involution $w$ on
$X_D^{(m)}$ is not an Atkin-Lehner involution and $\Aut
(X_D^{(m)})=\langle w\rangle \times W/\langle \om _m\rangle \simeq
\Z/2\Z \times \Z /2\Z $.

\end{prop}

{\bf Proof. } By the same arguments as in \cite{Ro}, Proposition
1, since $\Jac (X_D^{(m)})$ is isogenous to a product of simple
abelian varieties of real $\GL _2$-type, it follows that the group
of automorphisms of $X_D^{(m)}$ is abelian and only contains
involutions defined over $\Q $.

As it is checked from the genus formulas in \cite{Ogg1}, the $10$
cases in the above table are exactly those of Lemma \ref{g2} for
which the single Atkin-Lehner involution $u$ in $W/\langle \om
_m\rangle\subseteq \Aut (X_D^{(m)})$ is not the hyperelliptic
involution $w$, and hence is bielliptic. Since $\langle w,
u\rangle \subseteq \Aut (X_D^{(m)})\simeq (\Z /2\Z )^s$ for some
$s\geq 1$ and $g(X_D^{(m)})=2$, it follows that $s\leq 2$ and
hence $s=2$.

For the $29$ remaining cases not quoted in Proposition \ref{AL},
the single Atkin-Lehner involution in $W/\langle \om
_m\rangle\subseteq \Aut (X_D^{(m)})$ is the hyperelliptic
involution of the genus two curve $X_D^{(m)}$. We know from Kuhn \cite{Ku} that every quotient of a genus two curve $C/\Q $ by a
nonhyperelliptic involution defined over $\Q $ has a rational
point and thus is an elliptic curve over $\Q $.

Among these $29$ curves, it turns out by checking Cremona
tables, that their Jacobians are all simple over $\Q $ except for
$(D, m)= (57,3), (106,2)$ and $(118,2)$. Indeed, this follows from
the fact that these are the unique three cases such that there
exist two different isogeny classes of elliptic curves of
conductor $D$ and invariant by $\mu (m)\underline{\om }_m$. It is
then clear that these $26$ curves $X_D^{(m)}$ whose Jacobian is
simple over $\Q $ can not be bielliptic. As for the values $(D,
m)= (57,3), (106,2)$ and $(118,2)$ is concerned, there exactly two
isogeny classes of elliptic curves of level $D$ and invariant by
$\mu (m) \underline{\om }_m$. Namely, $57B, 57C$; $106A, 106C$;
$118B, 118C$, respectively. For each possible choice of elliptic
curves $E$ and $E'$ in these two isogeny classes, the abelian
surface $E\times E'$ contains no nontrivial $G_{\Q }$-stable
subgroups isomorphic to $(\Z /2 \Z )^2$. In other words, the
system (\ref{eq1}) admits no rational solution and therefore,
$X_D^{(m)}$ can not be a bielliptic curve. $\Box $

\begin{teor}\label{equations}
Equations for the curves in Proposition \ref{AL} and their
elliptic quotients are given in the following table:

$$
\begin{array}{|c|c|c|c|c|}
\hline D=p\cdot q & m & X_D^{(m)} & X_D^{(m)}/\langle \om
_q\rangle  &
X_D^{(m)}/\langle w\cdot \om _q\rangle   \\
\hline
91 = 7\cdot 13 & 91 & y^2=-x^6+19 x^4-3x^2+1  & 91B1  & 91A1 \\
\hline
123 = 3\cdot 41& 123 & y^2 =   - 9x^6 + 19x^4+ 5 x^2+1& 123A1   &  123B1 \\
\hline
141 = 3\cdot 47 & 141 &  y^2 = 27x^6 - 5x^4 - 7x^2 + 1 & 141A1  & 141D1  \\
\hline
142 = 2\cdot 71 & 2 & y^2 = -16x^6 - 87x^4 - 146x^2 - 71 & 142A1  & 142D2 \\
\hline
142 = 2\cdot 71 & 142 & y^2 = 16x^6 + 9x^4 - 10x^2 + 1 & 142A1  & 142B1  \\
\hline
155 = 5\cdot 31 & 155 &  y^2 = 25x^6 - 19x^4 + 11x^2 - 1 & 155A1  & 155C1  \\
\hline
158 = 2\cdot 79 & 158 & y^2 = -8x^6 + 9x^4 + 14x^2 + 1 & 158A1 & 158B1  \\
\hline
254= 2\cdot 127 & 254 &  y^2 = 8x^6 + 25x^4 - 18x^2 + 1  & 254A1 & 254C1 \\
\hline
326 = 2\cdot 163 & 326 &   y^2 =  x^6 + 10 x^4 - 63 x^2+4& 326B1  & 326A1 \\
\hline
446 = 2\cdot 223 & 446 &  y^2 = -16 x^6 - 7 x^4 + 38 x^2 + 1 & 446B1  & 446A1  \\
\hline

\end{array}
$$
Moreover, for all these equations the action of $w_q$  on them  is
$(x,y)\mapsto (-x,y)$.

\end{teor}

\noindent {\bf Proof.} We note that, in all $10$ cases, it follows
from Proposition \ref{AL} that $m=D$ except for the single case
$(D,m)=(142,2)$. Hence, the class of $\om _q$ in $\Aut
(X_D^{(m)})$ is the unique bielliptic involution of the curve that
is of Atkin-Lehner type. We have split the proof in five parts in
order to ease its reading.

\vskip 0.3cm {\em Step 1: Isogeny and isomorphism classes of the
elliptic quotients.}

We firstly determine the isogeny classes of the elliptic curves
$E=X_D^{(m)}/\langle \om _q\rangle$ and $E'=X_D^{(m)}/\langle
w\cdot \om _q\rangle $ of conductor $D$. There are exactly two normalized
newforms $f$, $f'$ for $\Gamma _0(D)$ with rational Fourier
coefficients whose signs of the eigenvalues for the action of the
Atkin-Lehner involutions are the following, depending on whether
$m=D$ or $m=p$:

$$
\begin{array}{|c|c|c|}
\hline
m=D & \underline{\om }_p & \underline{\om }_q \\
\hline f & - & - \\
\hline f' & + & + \\
\hline
 \end{array}
\qquad \qquad \qquad
\begin{array}{|c|c|c|}
\hline
m=p & \underline{\om }_p & \underline{\om }_q \\
\hline f & - & - \\
\hline f' & - & + \\
\hline
 \end{array}
$$

Then, the elliptic curves $E$ and $E'$ are isogenous to $A_f$ and
$A_f'$ over $\Q $, respectively. An examination of the $10$ cases
shows that, for $(D,m)\in \{ (141,141), (142,142),\\ (158, 158),
(326, 326), (446,446)\}$, the isomorphism classes of $E$ and $E'$
over $\Q $ are determined, because every isogeny class contains a
single isomorphism class. These are quoted in the table of the
statement. For the remaining cases, we have the following
possibilities:
$$
\begin{array}{|c|c|c|}
\hline
(D, m) & E & E' \\
\hline
(91,91) & B_1,B_2,B_3 & A_1 \\
\hline
(123, 123) & A_1, A_2 & B_1  \\
\hline
(142, 2) & A_1 & D_1, D_2 \\
\hline
(155, 155) &A_1, A_2 & C_1 \\
\hline
(254, 254) & A_1, A_2, A_3 & C_1  \\
\hline
 \end{array}
$$

\vskip 0.3cm {\em Step 2: Candidate equations.}

We now proceed to determine a finite set of candidate equations
for the $10$ curves $X_D^{(m)}$. We do so by applying Proposition
\ref{explicit} to every possible pair $(E, E')$ obtained in Step
$1$. For every pair, it turns out that system (\ref{eq1}) yields
 one rational solution. This means that, in the five cases
$(D, m)$ where there is a single possibility for $(E,E')$, we have
already determined an equation for the curve $X_D^{(m)}$, as
quoted in Theorem \ref{equations}. In the five remaining cases, we
obtain the following candidates:
$$
\begin{array}{|c|c|cccc|}
\hline
(D, m) & (E, E') &C: & a y^2&=&f(x) \\
\hline
(91,91) & (B_1,A_1) &  C_{91,1}:  & y^2&=&-x^6+19 x^4-3x^2+1 \\
    &    (B_2,A_1)  & C_{91,2}: & y^2& =&  91 x^6+ 43 x^4 + 9 x^2+1  \\
    &    (B_3,A_1)  & C_{91,3}: &  5y^2& =&  2401 x^6-403 x^4+3 x^4-1\\
\hline
(123, 123) & (A_1,B_1) &  C_{123,1}: &  y^2& =&- 9x^6 + 19x^4+ 5 x^2+1 \\
     &   (A_2,B_1)    &     C_{123,2}: &  y^2& =&   1681 x^6-419 x^4+35 x^2-1   \\
\hline
(142, 2) & (A_1,D_1) &  C_{142,1}: & y^2 &=& 8 x^6 + 33x^4 + 22x^2 + 1 \\
    &    (A_1,D_2)  & C_{142,2}: & y^2 &= &-16x^6 - 87x^4 - 146x^2 - 71 \\
\hline
(155, 155) &(A_1,C_1) &  C_{155,1}: & y^2 &=& 25x^6 - 19x^4 + 11x^2 - 1  \\
    &    (A_2,C_1)  & C_{155,2}: &  3 y^2 &=& 961 x^6 - 483x^4 - 45x^2 - 1 \\
\hline
(254, 254) & (A_1,C_1)    &  C_{254,1}: &  y^2& =& 8x^6 + 25x^4 - 18x^2 + 1  \\
        &   (A_2,C_1)    &  C_{254,2}: & y^2 & = & 127x^6 - 461x^4 - 51x^2 + 1  \\
       & (A_3,C_1)    &   C_{254,3}: & 71 y^2 & =& x^6 - 76888x^4 - 891x^2 + 2  \\
\hline
 \end{array}
$$
These equations have been computed such that the bielliptic
involution $\om _q$ acts as $(x, y)\mapsto (-x, y)$ and hence its
two fixed points are those whose $x$-coordinates are $0$. We
devote the rest of the proof to discard the wrong equations for
$X_D^{(m)}$ from the above in the table, by using suitable sieves.

\vskip 0.3cm {\em Step 3: The sieve of Heegner points fixed by the
involutions.}

Let us only consider in this part the four cases when $m=D$. From
Proposition \ref{AL}, $\Aut (X_D^{(D)}) = \langle \om _q, w\rangle
$. The two fixed points $P$, $Q$ of $\om _q$ acting on $X_D^{(D)}$
are the projection from $X_D$ onto $X_D^{(D)}$ of the four points
fixed by $\om _p$, when we regard it as an automorphism of $X_D$.
When $h(-4 p)=1$, it follows from the class field theory on
Heegner points (cf. \cite{JoPh}) that $P$, $Q\in X_D^{(D)}(\Q )$.
When $h(-4 p)=2$, it follows that $P$, $Q\in X_D^{(D)}(K)$, where
$K$ is a quadratic field such that $K(\sqrt {-p})$ is the Hilbert
class field of $\Q (\sqrt {-p})$. We thus obtain that, when $(D,
m)=(91,91)$, $(123,123)$ or $(254,254)$ the fixed points $P$, $Q$
lie on $X_D^{(D)}(\Q )$, whereas when $(D,m)=(155,155)$, $P$,
$Q\in X_{155}^{(155)}(K)$ for either $K=\Q (\sqrt {5})$ or $\Q
(\sqrt {-1})$.

On the other hand, on the bielliptic model $a y^2 = b x^6 + c x^4
+ d x^2 +e$, the coordinates of the two points fixed by $\om _q$
generate $\Q (\sqrt {e/a})$. This allows us to discard the
equations $C_{91,3}$, $C_{123,2}$, $C_{155,2}$ and $C_{254,3}$. We
have thus already determined an equation for $X_{123}^{(123)}$ and
$X_{155}^{(155)}$.

\vskip 0.3cm {\em Step 4: The real points sieve.}

Shimura proved in \cite{Sh2} that $X_D(\R )=\emptyset $. Later,
Ogg \cite{Ogg1} studied the question whether the Atkin-Lehner
quotients of Shimura curves admit real points. Namely, he proved
that $X_D^{(m)}(\R )\not =\emptyset $ if and only if
$(\frac{m}{p})\not =1$ for all $p\mid D$, $p\nmid m$.

Since $(\frac{2}{71})=1$, we deduce that $X_{142}^{(2)}(\R
)=\emptyset $ and hence $X_{142}^{(2)}\simeq C_{142,2}$ over $\Q
$.

\vskip 0.3cm {\em Step 5: The \v{C}erednik-Drinfeld sieve.}

Let us recall the theory of \v{C}erednik-Drinfeld on the bad
reduction of the Atkin-Lehner quotients $X_D^{(m)}$ at a fixed
prime $p\mid D$ (cf.\,\cite{Ba}, \cite{JoLi}, \cite{Kur},
\cite{Ogg2}).

Let $K_p$ denote the quadratic unramified extension of $\Q _p$ and
let $R_p$ be its ring of integers. Over $K_p$, the curves $X_D$
are generalized Mumford curves that admit a $p$-adic
uniformization by a Schottky group which is often
non-torsion-free. In the terminology of \cite{JoLi}, these are
called admissible curves.

Let $h(\delta ,\nu )$ denote the class number of a quaternion
Eichler order of level $\nu $ in a quaternion algebra of
discriminant $\delta $ over $\Q $. As shown in \cite{Kur}, $X_D$
admits a proper but often non regular integral model $\mathcal
{X}_D$ over $R_p$ whose special fibre $\widetilde {X}_D/\F _{p^2}$
is the union of $2 h(\frac {D}{p},1)$ irreducible rational
components meeting transversally at a total number of $h(\frac
{D}{p},p)$ points. The intersection points of the special fibre
are the only possible non regular points of $\mathcal{X}_D$ and
the only allowed multiplicities are $m=1,2$ and $3$. The reduction
type of $\mathcal{X}_D$ at $p$ is described by a weighted graph by
interpreting each component as a vertex, an intersection point $P$
between two components as an edge joining the two vertices and the
multiplicity $m$ of $P$ as the weight of the edge.

For every prime $q\leq 13$, it turns out that the dual graph of
$X_{p q}$ at $p$ consists of exactly two vertices joined by
$g(X_{p q})+1$ edges.

Moreover, the Atkin-Lehner involution $\om _{p q}$ lifts to
$\mathcal {X}_D$ and switches the two vertices and the quotient
graph consists of a single vertex with several loops of
multiplicity $1, 2$ or $3$ around it. In consequence, the special
fibre of $\mathcal {X}_D/\langle \om _{p q}\rangle $ has a single
and possibly singular irreducible component. After blowing up the
non regular closed points of $\mathcal {X}_D/\langle \om _{p
q}\rangle $ as in \cite[p.\ 288]{Kur}, we deduce that any two
irreducible components of the special fibre of the minimal regular
model of $X_D^{(m)}$ meet at most at two different intersection
points.

We can contrast this information with the explicit computation of
the reduction type of the equations in the above tables at the
primes $p\mid D$. This can be accomplished by means of Liu's
package {\em genus2reduction}, that computes the minimal regular
model of any curve of genus $2$ over $\Q $ over $\Z [\frac{1}{2}]
$.

The reduction type of $C_{91,1}$ and $C_{91,2}$ at $p=7$ are
$I_{\{ 1-1-0\} }$ and $I_{\{ 1-1-1\}}$, respectively. It
follows from \cite{NaUe}, the former is the symbol for a single
irreducible rational component with two nodes while the latter
corresponds to two rational components meeting at three points.
Hence, we discard $C_{91,2}$ and conclude that
$X_{91}^{(91)}\simeq C_{91,1}$. Similarly, the reduction type of
$C_{254,1}$ and $C_{254,2}$ at $p=127$ are $I_{\{1-1-0\}}$ and
$I_{\{1-1-1\}}$, respectively. This allows us to show that
$X_{254}^{(254)}\simeq C_{254,1}$. $\Box$

$\\ ${\em Acknowledgements.} The first author thanks the Number
Theory Group of the University of Nottingham for the warm
hospitality during the spring semester of 2003.

\end{document}